\documentclass[11pt]{scrartcl}

\usepackage[pdftex,x11names]{xcolor}
\usepackage[hidelinks]{hyperref}

\usepackage{chngcntr}
\usepackage{subfigure}
\usepackage[T1]{fontenc}
\usepackage[utf8]{inputenc}
\usepackage{lmodern}
\usepackage[english]{babel}

\usepackage{amsmath}
\usepackage{amssymb}
\usepackage{amsfonts}
\usepackage{amsthm}
\usepackage[centerdot]{mathtools}
\usepackage{dsfont}
\newcommand{\R}{\mathds{R}}
\newcommand{\Z}{\mathds{Z}}

\newcommand{\X}{\mathcal{X}}
\newcommand{\Y}{\mathcal{Y}}

\newcommand{\Ywsm}{\ensuremath{\Y_{S\lambda}}}
\newcommand{\Ybd}{\ensuremath{\Y_{S\partial}}}
\newcommand{\Ysnf}{\ensuremath{\Y_{SNF}}}

\usepackage{natbib}
\AtEndOfPackage{%
   \renewcommand\@openbib@code{%
      \advance\leftmargin\bibindent
      \itemindent -\bibindent
      \listparindent \itemindent
      \parsep \z@
      }%
   \renewcommand\newblock{\par}}%

\AtEndOfClass{\RequirePackage{natbib}%
\setlength{\bibhang}{\parindent}%

\def\@biblabel#1{#1.}%
\bibpunct{(}{)}{;}{a}{}{,}}

\usepackage{pifont}%

\DeclareMathOperator{\conv}{conv}

\usepackage{hyperref}
\usepackage[ruled,noend]{algorithm2e}
\SetKw{Continue}{continue}
\SetKw{Null}{null}
\SetKw{Break}{break}

\SetCommentSty{mycommfont}

\usepackage[textsize=tiny]{todonotes}
\usepackage[normalem]{ulem}
\usepackage{float}
\usepackage{authblk}

\usepackage{tikz}
\usepackage{tikz,tikz-3dplot}

\tikzset{
    cross/.pic = {
    \draw[rotate = 45] (-#1,0) -- (#1,0);
    \draw[rotate = 45] (0,-#1) -- (0, #1);
    }
}
\usetikzlibrary{shapes.misc}
\usetikzlibrary{shapes.geometric}
\usetikzlibrary{patterns}
\usetikzlibrary{patterns.meta,perspective}
\usetikzlibrary{graphs,quotes}
 \usepackage{pgfplots}
 \pgfplotsset{compat=1.16,
    every axis/.append style={
        axis lines=center,
        xlabel style={anchor=south west},
        ylabel style={anchor=south west},
        zlabel style={anchor=south west},
        tick align=outside,}
}
 \usepgfplotslibrary{patchplots}
\tdplotsetmaincoords{60}{120}

\usepackage{xpatch}

\newcommand{\weakly}{weakly }

\usepackage[]{cleveref} %nameinlink makes all blue

\newtheoremstyle{thm}% name
{15pt}% Space above
{5pt}% Space below 
{\itshape}% Body font
{}% Indent amount: Indent amount: empty = no indent, \parindent = normal paragraph indent
{\bfseries}% Theorem head font\part{title}
{}% Punctuation after theorem head
{0.5em}% Space after theorem head: Space after theorem head: { } = normal interword space; \newline = linebreak
{}% Theorem head spec (can be left empty, meaning `normal')
\theoremstyle{thm}% default

\newtheorem{thm}{Theorem}[section]
\newtheorem{lem}[thm]{Lemma}
\newtheorem{definition}[thm]{Definition}

\newtheorem{example}[thm]{Example}

\title{On Supportedness in Multi-Objective Combinatorial Optimization}

\author[1]{David K\"onen}
\author[1]{Michael Stiglmayr}
\affil[1]{%
	University of Wuppertal\\
	School of Mathematics and Natural Sciences\\
	Optimization Group\\
	Gaußstraße 20, 42103 Wuppertal, Germany\\ 
}
\affil[]{E-Mail:~\href{mailto:koenen@uni-wuppertal.de}{koenen@uni-wuppertal.de}, \href{mailto:stiglmayr@uni-wuppertal.de}{stiglmayr@uni-wuppertal.de} }

\date{}

\begin{document}

\maketitle

	\begin{abstract}\small

    This paper addresses an inconsistency in various definitions of supported non-dominated points within multi-objective combinatorial problems (\ref{eq:MOCO}). MOCO problems are known to contain supported and unsupported non-dominated points, with the latter typically outnumbering the former. Supported points are, in general, easier to determine, can serve as representations, and are used in two-phase methods to generate the entire non-dominated point set. Despite their importance, several different characterizations for supported efficient solutions (and supported non-dominated points) are used in the literature.  
    While these definitions are equivalent for multi-objective linear problems, they can yield different sets of supported non-dominated points for MOCO problems. We show by an example that these definitions are not equivalent for MOCO or general multi-objective optimization problems. Moreover, we analyze the structural and computational properties of the resulting sets of supported non-dominated points. These considerations motivate us to summarize equivalent definitions and characterizations for supported efficient solutions and to introduce a distinction between \emph{supported} and \emph{weakly supported} efficient solutions.

\end{abstract}

%%Graphical abstract
%\begin{graphicalabstract}
%\includegraphics{grabs}
%\end{graphicalabstract}

%%Research highlights
%\begin{highlights}
%\item Research highlight 1
%\item Research highlight 2
%\end{highlights}
	
	\par\vskip\baselineskip\noindent

\section{Introduction} 

\label{chapt:MODO}
Multi-objective optimization (MOO) is concerned with optimizing multiple conflicting objectives simultaneously. For a
comprehensive introduction to multi-objective optimization, we refer to~\cite{SteuerBook} and~\cite{ehrgott2005multicriteria}.

\begin{definition} A \emph{multi-objective optimization problem (MOO)} is defined as
	\begin{equation}\label{eq:MOP}
		\tag{MOO}
		\min_{x\in\X} f(x)=\left(f_1(x),\ldots,f_p(x)\right)^{\top},     
	\end{equation}
	where $f \colon \R^n \rightarrow \R^p$ is a vector-valued \emph{objective function} composed of $p\geq 2$ real-valued objective functions  $f_k \colon \R^n \rightarrow \R$ for $k \in \{ 1,\ldots,p\}$ and $\X \subseteq \R^n$ denotes the \emph{set of feasible solutions}. We call $\R^n$ the \emph{decision space} and $\R^p$ the objective space. The image of the feasible set  $$\Y \coloneqq f(\X) = \{ f(x) \colon x\in \X\}\subseteq \R^p$$ is called the set of \emph{feasible outcome points} in objective space. 
\end{definition}

The non-negative orthant of $\R^p$ is denoted by $\mathbb{R}_{\geqq}^p:=\left\{x \in \mathbb{R}^p: x \geqq 0\right\}$ and analogously its interior $\R_{>}^p$ and $\R_{\geq}^p := \left\{x \in \R^p: x \geqq 0, x \neq 0 \right\}$. In our notation, we also use the Minkowski sum and the Minkowski product of two sets $A, B \subseteq \R^{p}$, which are defined as $A+B\coloneqq \{a+b \colon a \in A, b \in B\}$ and $A \cdot B\coloneqq \{a \cdot b \colon $ $a \in A, b \in B\}$, respectively. For a closed set $C$, we define the boundary of $C$ by $\partial C$. 

The objective functions are typically assumed to be conflicting, implying the absence of an \emph{ideal} solution that minimizes all objectives simultaneously.  Instead, we rely on the \emph{Pareto concept of optimality}, based on the component-wise order in $\R^p$.

	\begin{definition}
		\label{def:component}
		Let $y^1, y^2 \in \mathbb{R}^p$. We write: 
		\begin{itemize} 
			\item $y^1 \leqq y^2$ if $y_i^1 \leq y^2_i$ for $i=1, \ldots, p$, 
			\item $y^1 \le y^2$ if $y^1 \leqq y^2$ but $y^1 \ne y^2$ and 
			\item $y^1 < y^2$ if $y_i^1 < y^2_i$ for $i=1, \ldots, p$. 
		\end{itemize}
\end{definition}

\begin{definition}
	A feasible solution $x^{*} \in \X$ is called \emph{efficient} if there is no $x \in \X$ such that $f(x) \leq f\left(x^{*}\right)$. The image $f(x^*)$ is then called a \emph{non-dominated point}. 
	The set of efficient solutions is denoted by $\X_E\subseteq \X$ and the set of non-dominated points by $\mathcal{Y}_N\subseteq \Y$.  
	
	A feasible solution $\hat{x}\in X$ is called \emph{weakly efficient} if no $x \in \X$ exists with $f(x) < f(\hat{x})$. The image $f(\hat{x})$ is then called a \emph{weakly non-dominated point}. 
	Furthermore, a feasible solution $x^* \in \mathcal{X}$ is called \emph{properly efficient} (in Geoffrion's sense \citep{GEOFFRION1968618}), if it is efficient and if there is a real number $M>0$ such that for all $k\in\{1,\ldots,p\}$ and $x \in \mathcal{X}$ satisfying $f_k(x)<f_k({x}^*)$ there exists an index $j$ such that $f_j({x}^*)<f_j(x)$ such that
	$$
	\frac{f_k({x}^*)-f_k(x)}{f_j(x)-f_j({x}^*)} \leq M .
	$$
	The image ${y}^*=f({x}^*)$ is called properly non-dominated.
\end{definition}

In the following, we assume that the multi-objective optimization problem has at least one non-dominated point, i.e., $\Y_N\neq \varnothing$. 
The polyhedron $\mathcal{Y}^{\geqq}:=\conv(\mathcal{Y}_N+\R^p_{\geqq})$ is called the \emph{upper image} of \(\mathcal{Y}\).

Scalarization approaches are the most common solution techniques in multi-objective optimization, which rely on replacing the multi-objective problem by a sequence of single-criteria optimization problems.
In the \emph{weighted-sum scalarization} a convex combination (with weights \(\lambda_i\)) of the objective functions is optimized over the feasible set \(\mathcal{X}\). 
%Let $\|x\|_{1}$ denote the $1$-norm of $x \in \R^{d}$, i.\,e.,  $\|x\|_{1}\coloneqq \sum_{i=1}^{d}\left|x_{i}\right|$.
Then, the set of \emph{normalized weighting vectors} (or \emph{weight space}) is defined as the set $\Lambda_p=\{\lambda \in \R^{p}_{>} \colon \sum_{i=1}^p \lambda_i = 1 \}$ or $\Lambda_p^0=\{\lambda \in \R^{p}_{\geq} \colon \sum_{i=1}^p \lambda_i = 1 \}$ if weights equal to zero are included.

\begin{definition}
	The \emph{weighted-sum scalarization} of \eqref{eq:MOP} with weighting vector $\lambda \in \Lambda_p$ or $\lambda \in \Lambda_p^0$ is defined as the parametric program 
	\begin{equation}
		\tag{$P_{\lambda}$}\label{eq:wsm}
		P_{\lambda}\coloneqq  \min_{x\in \X}\;  \lambda^\top f(x) 
	\end{equation}
\end{definition}

\begin{thm}[\citealt{ehrgott2005multicriteria}]\label{thm:ehrgott}
	If $\lambda\in \Lambda_{p}^0$, every optimal solution of \ref{eq:wsm} is a weakly efficient solution of \eqref{eq:MOP}. Moreover, every optimal solution of \ref{eq:wsm} problem is properly efficient for  \eqref{eq:MOP}, if \(\lambda\in\Lambda_{p}\).
\end{thm}

A specific, well-studied class of \ref{eq:MOP} are \emph{multi-objective linear programs (MOLP)} and \emph{multi-objective integer linear programs (MOILP)}.

\begin{definition} A \emph{multi-objective linear optimization problem (MOLP)} is defined by
	\begin{equation}\label{eq:MOLP}
		\tag{MOLP}
		\min_{x\in\X} y(x)=\left(y_1(x),\ldots,y_p(x)\right)^{\top} = Cx,     
	\end{equation}
	with $\X\coloneqq \{   x  \in \R^n  \colon Ax \leqq b, x\geqq 0 \}$. Thereby, the rows $c^k$ ($k=1,\ldots,p$) of the cost matrix $C \in \mathbb{R}^{p \times n}$ contain the coefficients of the $p$
	linear objective functions $y_k(x)=c^k\, x$ for $ k\in \{1, \ldots, p\}$ and the matrix $A \in \mathbb{R}^{m\times n}, b \in \mathbb{R}^m$ formulates the $m$ linear constraints. 
\end{definition}

\begin{definition}
	If all decision variables are further restricted to the set of integers, i.e., $x \in \Z$, we obtain
	a multi-objective integer linear problem (MOILP).
	\begin{equation}\label{eq:MOILP}
		\tag{MOILP}
		\min_{x\in\X} y(x)=\left(y_1(x),\ldots,y_p(x)\right)^{\top} = Cx,     
	\end{equation}
	with $\X\coloneqq \{   x  \in \Z^n  \colon Ax \leqq b, x\geqq 0 \}$ and all other elements as described in \eqref{eq:MOLP}. 
\end{definition}

\begin{thm}[\citealt{Iser74m}]\label{thm:iser}
	Considering \ref{eq:MOLP} problems, for every efficient solution $\bar{x}\in\mathcal{X}_E$ there exists a $\lambda \in \Lambda_p$ such that $\bar{x}$ is a optimal solution of the weighted sum scalarization \eqref{eq:wsm}.
\end{thm}
Hence, for a~\ref{eq:MOLP}, the efficient set $\X_E$ and the union over all sets of optimal solutions of \ref{eq:wsm} for $\lambda\in \Lambda_p$ are identical. Thus, in multi-objective linear programs (MOLP) every non-dominated point lies on the non-dominated frontier, which can be defined in the following way. 

\begin{definition}[\cite{Hamacher2007}]\label{def:non-dominated-frontier}
	The \emph{non-dominated frontier} is the set 
	$$\Bigl\{y \in \conv(\Y_N) \colon \conv(\Y_N) \cap (y -\R^p_{\geqq}) = \{y\} \Bigr\}.$$
\end{definition}
Note that Hamacher uses the wording efficient
frontier instead of the non-dominated frontier. The non-dominated frontier equals the non-dominated set of $\conv (\Y_N)$ and can be characterized as the set of the \emph{maximal non-dominated faces} of the upper image $\Y^{\geqq}$~\citep{ehrgott2005multicriteria}. Hereby, a \emph{non-dominated face} $F\subseteq \Y^{\geqq}$ is a face of the upper image \(\Y^{\geqq}\) %convex hull of the feasible outcome set 
such that every point on $F$ is  non-dominated with respect to $\Y^{\geqq}$. Similarly, $F\subseteq \Y^{\geqq}$ is called a \emph{weakly non-dominated face} if all its points are weakly non-dominated. A face $F$ is called \emph{maximally non-dominated} if there is no other non-dominated face $G$ of \(\Y^{\geqq}\) such that $F$ is a proper subset of $G$. This would imply that the dimension of $G$ is greater than the dimension of $F$.  The preimage $F_\X$ of a maximally non-dominated face $F_\Y$ of \(\Y^{\geqq}\), i.e., all solutions whose image lies in $F_\Y$, is denoted as the \emph{maximally efficient face}. Note that $F_\X$ might not be a face of the feasible set, the polyhedron \(\X\), since multiple feasible solutions (located on different faces of $\X$) could map to the same non-dominated point $y\in\Y$. 

\begin{figure}[h]
	\centering
	\tdplotsetmaincoords{120}{110}
	\begin{tikzpicture}[tdplot_main_coords,line join=round,
		declare function={a=4; b=2; Tx=2; Ty=3; Tz=6;}, scale=.7]
		
		% coordinate axes
		\draw (0,0,0) -- (8,0,0);
		\draw (0,0,0) -- (0,8,0);
		\draw (0,0,0) -- (0,0,8); 
		\draw[->] (8,0,0) -- (10,0,0) node[anchor=south]{$y_3$};
		\draw[->] (0,8,0) -- (0,10,0) node[anchor=north west]{$y_1$};
		\draw[->] (0,0,8) -- (0,0,10) node[anchor=south]{$y_2$};
		
		% hidden edges (translated)
		\draw[thick]  (Tx, Ty+ b, Tz- a) -- (Tx, Ty+ a, Tz- a);
		\draw[thick]  (Tx, Ty, Tz) -- (Tx, Ty, Tz- a+ b); 
		
		% visible edges
		\draw[thick] (Tx+ a, Ty, Tz) -- (Tx+ a, Ty, Tz- a);
		\draw[dashed] (Tx+ a, Ty, Tz- a) -- (Tx+ a, Ty+ a, Tz- a);
		\draw[thick] (Tx+ b, Ty, Tz- a) -- (Tx+ a, Ty, Tz- a);
		\draw[thick]  (Tx+ a, Ty+ a, Tz) -- (Tx+ a, Ty, Tz) -- (Tx, Ty, Tz) 
		-- (Tx, Ty+ a, Tz) -- (Tx+ a, Ty+ a, Tz);
		\draw[dashed]  (Tx, Ty+ a, Tz- a) 
		-- (Tx+ a, Ty+ a, Tz- a) -- (Tx+ a, Ty+ a, Tz);
		
		\draw[thick] (Tx, Ty+ a, Tz) -- (Tx, Ty+ a, Tz- a);      
		% cut face 
		\draw[thick,pattern={Lines[angle=45,distance={4.5pt}]}] 
		(Tx, Ty, Tz- a+ b) -- (Tx+ b, Ty, Tz- a) -- (Tx, Ty+ b, Tz- a) -- cycle;

		\fill[cyan!40,opacity=0.6]
		(Tx, Ty, Tz) -- (Tx, Ty+a, Tz) -- (Tx, Ty+a, Tz-a) -- (Tx, Ty+b, Tz-a) -- (Tx, Ty, Tz-a+b) --cycle;

		\fill[cyan!40,opacity=0.6]
		(Tx, Ty, Tz) -- (Tx+a, Ty, Tz) -- (Tx+a, Ty, Tz-a) --  (Tx+b, Ty, Tz-a)  -- (Tx, Ty, Tz-a+b)  -- cycle ; 
		
		%Boundaries
		\draw[thick,cyan!60]
		(Tx, Ty, Tz) -- (Tx, Ty+a, Tz) -- (Tx, Ty+a, Tz-a) -- (Tx, Ty+b, Tz-a) -- (Tx, Ty, Tz-a+b) -- cycle;
		
		\draw[thick,cyan!60]
		(Tx, Ty, Tz) -- (Tx+a, Ty, Tz) -- (Tx+a, Ty, Tz-a) -- (Tx+b, Ty, Tz-a) -- (Tx, Ty, Tz-a+b) -- cycle;

		\fill[cyan!40,opacity=0.6]
		(Tx+b, Ty, Tz-a) -- (Tx+a, Ty, Tz-a) --   (Tx+a, Ty+a, Tz-a) --  (Tx, Ty+a, Tz-a) --  (Tx, Ty+b, Tz-a) -- cycle; 
		
		% \draw[thick,cyan!60]
		%   (Tx+b, Ty, Tz-a) -- (Tx+a, Ty, Tz-a) --   (Tx+a, Ty+a, Tz-a) --  (Tx, Ty+a, Tz-a) --  (Tx, Ty+b, Tz-a) -- cycle; 
		
		\fill[red!40,opacity=0.6] 
		(Tx, Ty, Tz- a+ b) -- (Tx+ b, Ty, Tz- a) -- (Tx, Ty+ b, Tz- a) -- cycle;
		\draw[thick,red!70] 
		(Tx, Ty, Tz- a+ b) -- (Tx+ b, Ty, Tz- a) -- (Tx, Ty+ b, Tz- a) -- cycle;
		
		% edges of the cut face (non-dominated faces)
		\draw[line width=1.2pt,red!60] (Tx, Ty, Tz- a+ b) -- (Tx+ b, Ty, Tz- a);
		\draw[line width=1.2pt,red!60] (Tx+ b, Ty, Tz- a) -- (Tx, Ty+ b, Tz- a);
		\draw[line width=1.2pt,red!60] (Tx, Ty+ b, Tz- a) -- (Tx, Ty, Tz- a+ b);

	\end{tikzpicture}
	\caption{ Illustration of dominance relations in faces: the red triangular face represents a \emph{maximal non-dominated face}; its red boundary edges are \emph{non-dominated faces}; and the light blue faces are \emph{weakly non-dominated}.}
	\label{fig:faces}
\end{figure}

 Figure~\ref{fig:faces} illustrates the geometric structure of non-dominated and weakly non-dominated faces for a continuous multi-objective minimization problem, where the feasible region consists of a cube with a cut at its lower-left front corner. The red triangular region corresponds to a \emph{maximal non-dominated face}, meaning no point in this region is dominated by another feasible point, and this face cannot be extended to a larger non-dominated face within the feasible set. The bold red edges of this triangle represent lower-dimensional \emph{non-dominated faces}, while the light-blue shaded cube faces indicate \emph{weakly non-dominated faces}, which are optimal in some but not all directions.

In contrast to multi-objective linear programs (MOLP), where every non-dominated point lies on the non-dominated frontier, and thus on the convex hull of the upper image, MOILPs or \ref{eq:MOCO} problems may contain non-dominated points that are located in the interior of the upper image, i.e., in $\operatorname{int}(\conv(\Y_N + \R^p_{\geqq}))$ and cannot be generated through any weighted-sum scalarization.

\begin{definition} Formally a \emph{multi-objective combinatorial optimization} (MOCO) problem is defined as
	\begin{equation}\label{eq:MOCO}
		\tag{MOCO}
		\min_{x\in\X}  y(x)=\left(y_1(x),\ldots,y_p(x)\right)^{\top} = C\, x,     
	\end{equation}
	where $\X \coloneqq \{ x  \in  \{0,1\}^n \colon Ax = b \}$, 
	with the cost matrix $C \in \mathbb{Z}^{p \times n}$ containing the rows $c^k$ of coefficients of $p$
	linear objective functions $y_k(x)=c^k\, x$ for  $k \in \{1, \ldots,p \}$ and $A \in \mathbb{Z}^{m\times n}, b \in
	\mathbb{Z}^m$ describing the $m$ constraints. The constraints define combinatorial structures such as paths, trees, or cycles 
	in a network or partitions of a set. 
\end{definition}

This definition includes various problems such as multi-objective spanning trees, shortest paths, knapsack and assignment problems.

Two major challenges are associated with MOCO problems. Firstly, they are theoretically more challenging than their single-objective counterparts, falling into the class of computationally intractable problems~\citep{ehrgott2005multicriteria}. Secondly, in MOCO problems, non-dominated points located in the interior of the upper image can occur,  which typically outnumber the set of solutions on the non-dominated frontier (see, e.g.,~\citealt{Visee1998}). Hence, a distinction between different classes of efficient solutions is required.
An intuitive definition of supportedness is the following.

\begin{definition}\label{def:supp1} 
	\begin{enumerate}
		\item An efficient solution is called \emph{supported efficient solution} if it is an optimal solution to the weighted sum scalarization \ref{eq:wsm} for  $\lambda\in \Lambda_{p}$, i.e., an optimal solution to a single objective weighted-sum problem where the weights are strictly positive. Its image is called \emph{supported non-dominated point}; we use the notation $\X_{S}$ and $\Y_{S}$ for the supported efficient solution set and supported non-dominated point set, respectively. For discrete problems, supported non-dominated points are located on the \emph{non-dominated frontier}, i.e., located on the union of the \emph{maximal non-dominated faces} of the upper image.   
		
		\item \emph{Extreme supported efficient solutions}  are those solutions whose image lies on the vertex set of the upper image. Their image is called an \emph{extreme supported non-dominated point}. 
		%We use the notation $\Y_{EE}$ for the set of extreme non-dominated points and $\X_{E}$ for the set of extreme efficient solutions.
		
		\item \emph{Unsupported efficient solutions} are efficient solutions that are not optimal solutions of  \ref{eq:wsm} for any  $\lambda\in \Lambda_{p}^0$. Unsupported non-dominated points lie in the interior of the upper image. 
	\end{enumerate}
\end{definition}

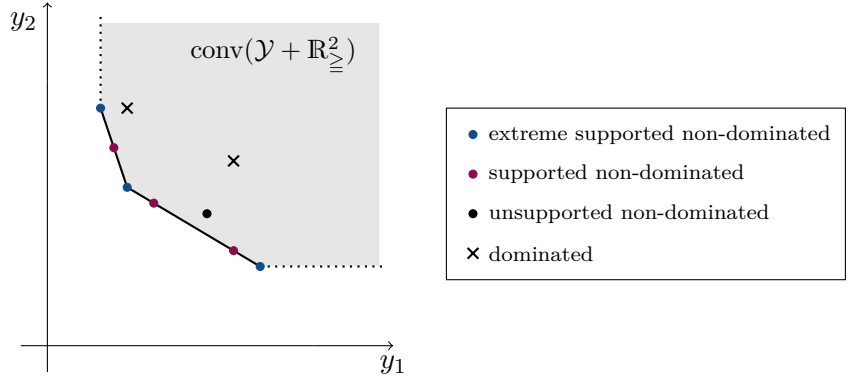
\begin{figure}
	\centering
	\begin{tikzpicture}[]
		\tdplotsetmaincoords{0}{0}
		\tikzstyle{vertex}=[diamond,fill= black,draw=black,minimum size=3.5pt,inner sep=0]
		\tikzstyle{vertex2}=[rectangle,fill=DeepPink4,draw=DeepPink4,minimum size=3pt,inner sep = 0]
		\tikzstyle{vertex3}=[circle,fill=DodgerBlue4,draw=DodgerBlue4,minimum size=3.5pt,inner sep = 0]
		\tikzstyle{N_point}=[draw, cross out,scale=.5,thick]
		
		\begin{scope}[tdplot_main_coords, scale=.35]
			
			%Shaded Area:
			\draw [draw= white, fill= black!10  ,fill opacity=1 ]  
			(2,12.25) -- (2,9) -- (3,6)  -- (8,3) -- (12.5,3) -- (12.5,12.25)-- cycle;
			
			\draw[] (-1,0) -- (10,0);
			\draw[] (0,-1) -- (0,10);
			\draw[->] (10,0) node[anchor=north]{} -- (13,0) node[anchor=north]{{$y_1$}};
			\draw[->] (0,10) -- (0,13) node[anchor=north east]{{$y_2$}};
			
			\draw[thick] (2,9) -- (3,6) -- (8,3);
			\draw[dotted,thick] (2,9) -- (2,12.5);
			\draw[dotted,thick] (8,3) -- (12.75,3);
			
			\node[] at ( 8.5 ,11 ){\small $ \conv(\Y_N+\R^2_{\geqq}) $};
			
			\node[vertex3, label={[label distance=-5pt]above right:{}}] at (2,9) {};
			\node[vertex3, label={[label distance=-5pt]above right:{}}] at (3,6) {};
			\node[N_point, label={[label distance=-5pt]above right:{}}] at (6,5) {};
			\node[vertex3, label={[label distance=-5pt]above right:{}}] at (8,3) {};
			\node[vertex2, label={[label distance=-5pt]above right:{}}] at (4,5.4) {};
			\node[vertex2, label={[label distance=-5pt]above right:{}}] at (7,3.6) {};
			\node[vertex2, label={[label distance=-5pt]above right:{}}] at (2.5,7.5) {};
			\node[vertex, label={[label distance=-5pt]above right:{}}] at (3,9) {};
			\node[vertex, label={[label distance=-5pt]above right:{}}] at (7,7) {};
			
			\draw (15,2.5) -- (15,9) -- (30,9) -- (30,2.5) -- cycle;
			\node[vertex3] at (16,8) {};
			\draw (16.5,8) node [anchor= west][inner sep=0.75pt]   [align=left] {\scriptsize extreme supported non-dominated};
			\node[vertex2] at (16,6.5) {};
			\draw (16.5,6.5) node [anchor= west][inner sep=0.75pt]   [align=left] {\scriptsize supported non-dominated};
			\node[N_point] at (16,5) {};
			\draw (16.5,5) node [anchor= west][inner sep=0.75pt]   [align=left] {\scriptsize unsupported non-dominated};
			\node[vertex] at (16,3.5) {};
			\draw (16.5,3.5) node [anchor= west][inner sep=0.75pt]   [align=left] {\scriptsize dominated};
		\end{scope}
		
	\end{tikzpicture}
	\caption{Illustration of the upper image $\mathcal{Y}^{\geqq}=\conv(\mathcal{Y})+ \R^2_{\geqq}$ and the different solution types.}\label{Image:DiffSolutions}\end{figure}

\Cref{Image:DiffSolutions} illustrates supported extreme, supported, and unsupported non-dominated points, as well as the upper image in the bi-objective case.

As unsupported non-dominated points cannot be computed through a weighted sum scalarization, different scalarization techniques are required that are often computationally more expensive.

The determination of supported non-dominated points has gained attention for several reasons. First, these points are generally easier to determine than the unsupported non-dominated points and, moreover, the unsupported non-dominated points in MOCO typically outnumber the supported ones, as observed in Studies as~\cite{Visee1998}.  

Second, they can serve as a foundation for the second phase of \emph{two-phase methods}, which aims to generate the entire non-dominated point set using information derived from the supported non-dominated points (\citealt{pasternak1972bicriterion,Visee1998, Hamacher2007, PryEtAl2008,PRZYBYLSKI2010149,EUSEBIO20092554, dai2018two},\ldots). Recently, a computational study \cite{Serpil2024} showed that the set of supported non-dominated points can be used as a high-quality representation in multi-objective discrete optimization problems, focusing specifically on binary knapsack and assignment problems.

Several studies focus on identifying or analyzing the supported or extreme supported point set across various combinatorial problems, including bi- and multi-objective integer network flow problems~\citep{EUSEBIO200968,Church2015,raith17,konen2023outputpolynomial,konen2023outputsensitive}, bi- and multi-objective minimum spanning trees~\citep{sourd2006multi,silva07note,Correia2021}, bi- and multi-objective shortest path problems~\citep{EDWIN99,SEDENONODA15,church14}, bi- and multi-objective combinatorial unconstrained problems~\citep{boeklerthesis,Schulze2019}, bi- and multi-objective knapsack problems~\citep{Visee1998,argyris2011identifying,schulze2017new}, and bi- and multi-objective assignment problems~\citep{tuyttens2000performance,gandibleux2003use,PRZYBYLSKI2010149}. Additional work on general MOCO problems can be found in~\cite{gandibleux01,jesus2015implicit,Serpil2024}, among others.

Beyond MOCO problems, there is research on the identification of supported solutions in multi-objective mixed integer problems~\citep{oezpeynirci10exact,ozlen2019,Bokler2024} or addressing supportedness in non-convex problems~\citep{dinh1995,liefooghe2014hybrid,liefooghe2015experiments}. 
Summarizing, supported points are often easier to determine, can serve as high-quality representations, and can be used in two-phase methods to generate the entire non-dominated point set.

However, despite their significance, several alternative definitions (as in \Cref{def:supp1}) and characterizations for supportedness are  used. These different definitions differ in the literature and sometimes even within a single publication. 
The following chapter proves that these definitions are not equivalent in the case of discrete or more general optimization problems. As a result, they  generate different sets of supported efficient solutions and, consequently, different supported non-dominated point sets with varying properties.

%Supported non-dominated points for MOCO or MOILP problems are often characterized as non-dominated points that lie on the boundary of the upper image and that they only lie on the \emph{non-dominated frontier}, %which can be characterized as the union of all \emph{maximally non-dominated faces}, 
%whereas unsupported solutions are non-dominated points that lie in the interior of the upper image~\citep{EUSEBIO200968,Przy10}.  
%
%%Hereby, a face of the upper image is called a \emph{non-dominated face} if all its points are non-dominated. A face $F$ is called \emph{maximally non-dominated} if there is no other non-dominated face $G$ such that $F\subset G$, this would imply that the dimension of the face $G$ is greater than the dimension of the face $F$. \david{faces weglassen? ich habe 2x faces eingefügt!}
%
%However, depending on the definition, unsupported non-dominated points may exist that lie on the boundary of the upper image or supported non-dominated points that lie on \emph{weakly non-dominated faces} and thus do not lie on the non-dominated frontier. In this case, they cannot be obtained as optimal solutions of a weighted sum problem with weights strictly greater than zero.

This motivates us to distinguish between \emph{supported} efficient solutions and \emph{weakly supported} efficient solutions. 
An efficient solution is denoted as weakly supported if it is an optimal solution of a weighted sum scalarization with non-negative weights. In contrast to the definition of supportedness, weakly supportedness allows single weights to have a value of zero.
 Note that the terminology weakly supported efficient solutions, despite the fact that these solutions are still efficient, is motivated by geometric considerations: any weakly supported non-dominated point that is not a supported non-dominated point lies on a weakly non-dominated face, but not on a non-dominated face. In contrast, supported points necessarily lie on non-dominated faces.

% Then the following characterizations for MOILP problems hold: While weakly supported non-dominated points lie on the boundary of the upper image, supported non-dominated points lie only on the non-dominated frontier, i.\,e., only on maximally non-dominated faces. The unsupported non-dominated points lie in the interior of the upper image. 
%We will present an example where the set of supported efficient solutions is a proper subset of all weakly efficient solutions in~\Cref{chapt:supp_vs_strictly}. The clear distinction between the sets of supported non-dominated and weakly supported non-dominated solutions is also necessary as the corresponding problems may differ in their output time complexity. In particular, in the case of the minimum cost flow problem, it can be shown that supported efficient solutions can be determined in \emph{output-polynomial time}, whereas this is not the case for the \emph{weakly supported solutions} unless $\mathbf{P} = \mathbf{NP}$~\citep{konen2023outputpolynomial}. 

This work focuses on the theoretical understanding of supportedness in MOCO problems. In particular, we analyze the structural properties of supported and weakly supported non-dominated points and highlight important distinctions through counterexamples. However, further investigation into the supportedness of general MOO problems remains sparse and presents an important area for future research. A recent study in~\citep{ChlumskyHarttmann2025} provided a first approach toward a categorization of supportedness definitions for general MOO problems.

The remainder of the paper is structured as follows. In~\Cref{chapt:supp_vs_strictly} the different definitions and characterizations of supportedness found in the literature are presented. A counterexample is given that shows that these definitions are not equivalent for MOILPs and MOCO problems. Additionally, it discusses the distinction between supported and weakly supported non-dominated points. Finally,~\Cref{chapt:concl} concludes the findings and outlines directions for future research.

\section{Supported and Weakly Supported Non-dominated Points in MOCO problems}\label{chapt:supp_vs_strictly}

Several definitions and characterizations despite~\Cref{def:supp} of supported non-dominated points exist in the literature and sometimes are even inconsistent within a single publication.

\begin{definition}[Conflicting Definitions of Supportedness]\label{def:supp}
	A point $y'\in \Y$ is called supported 
	\begin{itemize}
		\item[(1)]\label{supp2} if $y'$ is non-dominated  and $y'$ lies on the boundary of the upper image $\mathcal{Y}^{\geqq}\coloneqq \conv(\mathcal{Y}_N+\R^p_{\geqq })$, i.e., $y'\in \Y_N\cap \partial \Y^{\geqq}$~\citep{Eusebio09,EUSEBIO200968,liefooghe2014hybrid,liefooghe2015experiments,Correia2021}. The corresponding set of supported non-dominated points is denoted by $\Y_{S\partial}$.
		\item[(2)]\label{supp3} if $y'$ is non-dominated and $y'$ is located on the \emph{non-dominated frontier} defined as the set $\{y \in \conv(\mathcal{Y}_N) \colon \conv(\mathcal{Y}_N) \cap (y -\R^p_{\geqq}) = \{y\}\}$~\citep{Hamacher2007}, the set supported non-dominated points according to this definition is denoted by $\Y_{SNF}$.
		\item[(3)]\label{supp4} if $y'$ is the image of a supported efficient solution, which are those efficient solutions that can be obtained as optimal solutions of a weighted sum scalarization with weights strictly greater zero, i.e. $y$ is an image of a solution $x \in \arg \min P_\lambda$ for a $\lambda\in \Lambda_p$~\citep{Visee1998,ehrgott2005multicriteria,raith09,Przy10,argyris2011identifying,raith17}. The set of supported non-dominated points obtained with respect to this definition is denoted by $\Y_{S\lambda}.$
	\end{itemize}
	%Let $\Y_{S_i}$ be the set of supported non-dominated points obtained regarding the $i$-th definition with $i\in\{1\ldots,4\}$.
\end{definition}

Note, that some publications adopt the definition related to $\Y_{S\lambda}$ but allowing $\lambda\in \Lambda_p^0$, provided that $y$ is non-dominated, see, e.g.~\cite{gandibleux01,liefooghe2015experiments}. Regarding definition $ \Y_{SNF}$,~\cite{Hamacher2007} uses the wording \emph{efficient frontier} instead of the non-dominated frontier. There also exist definitions based on convex combinations of the non-dominated points as given in~\cite{oezpeynirci10exact}.

\cite{Iser74m} showed that for a~\ref{eq:MOLP}, the set of non-dominated efficient solutions coincides with the set of optimal solutions of the weighted sum method with weights strictly greater than zero, i.e., in ~\ref{eq:MOLP}, the weighted sum scalarization achieves completeness when varying $\lambda$ within $\Lambda_p$. Consequently, for a ~\ref{eq:MOLP}, it holds $\Y_{S\partial} =  \Y_{SNF} =  \Y_{S\lambda}  $.

\begin{thm}
	For a~\ref{eq:MOLP} it holds $\Y_{S\partial} =  \Y_{SNF} =  \Y_{S\lambda}  $.
\end{thm}

\begin{proof}
	According to \Cref{thm:iser}, any non-dominated point $y'\in \Y$ in~\ref{eq:MOLP} can be obtained as the image of an optimal solution $x\in \arg\min P_{\lambda} $ with $\lambda\in \Lambda_p$, and with \Cref{thm:ehrgott} we obtain $\Y_{S\lambda} = \Y_N$. Any point $\hat{y}\in \Y_{S\partial}$ or $\hat{y}\in \Y_{SNF}$  is non-dominated per definition and therefore $\hat{y}\in \Y_{S\lambda}.$ It follows $$\Y_{S\lambda} \supseteq \Y_{S\partial} \text{ and } \Y_{S\lambda} \supseteq  \Y_{SNF}.$$
	
	Note that for an \ref{eq:MOLP}, $\Y$ is polyhedral and thus $\Y = \conv(\Y)$. Consider $\hat{y}\in \Y_{S\lambda}$, i.e.,  $\hat{y}$ is an image of an optimal solution $x' \in \arg\min P_{\lambda}$ with $\lambda \in \Lambda_p$. Since $\Y$ is polyhedral, $\hat{y}$ must lie in a face of $ \conv (\Y)$ \citep{nem99}. Hence, $\hat{y} \in \partial \conv (\Y)  $ and from~\Cref{thm:ehrgott} it follows $\hat{y}\in \Y_N$. Thus, it follows $\hat{y}\in \partial \conv(\Y_N)\cap \Y_N$, and hence $\hat{y}\in \Y_{\partial}$. Furthermore, since $\conv (\Y_N)$ is a nonempty convex set, it holds  that $\hat{y}$ is a minimal element of $\conv (\Y_N)$ induced by the cone $-\R^p_{\geqq}$ \citep{boyd2004convex} and we can conclude that $((\hat{y}-\R^p_{\geqq})\backslash{\hat{y}})\cap \conv(\Y_N) = \varnothing$. Consequently, it follows that $\hat{y}\in \Y_{SNF}$. Together it holds,
	$$
	\Y_{S\lambda} \subseteq \Y_{S\partial} \text{ and } \Y_{S\lambda} \subseteq \Y_{SNF} . 
	$$
	
	This concludes $\Y_{S\partial} =  \Y_{SNF}  = \Y_{S\lambda} $.
\end{proof}

This equivalence does not hold for the discrete or the non-linear case.
In the literature, supported non-dominated points for~\ref{eq:MOCO} or~\ref{eq:MOILP} 
are often characterized as non-dominated points that lie on the boundary of the upper image
and that they only lie on the maximally non-dominated faces, i.e., the non-dominated frontier. In contrast, the unsupported are characterized as non-dominated points that lie in the interior of the upper image, e.g.,~\citep{EUSEBIO200968,Przy10,Correia2021}.

However, depending on the definition, unsupported non-dominated points may exist that lie on the boundary of the upper image or supported non-dominated points that lie on \emph{weakly non-dominated faces} and thus do not lie on the non-dominated frontier. In this case, they cannot be obtained as optimal solutions of a weighted sum problem with weights strictly greater than zero.
Suppose a supported non-dominated point is defined according to \Cref{def:supp} (2) and (3), i.e., considering the sets $\Y_{S\lambda}$ and $\Y_{SNF}$. In that case, unsupported points may exist on the boundary of the upper image. Conversely, if a supported non-dominated point is defined according to (1), i.e., points contained in \Ybd, there may exist supported non-dominated points on the boundary of the upper image which are not lying on the non-dominated frontier, i.\,e., which lie on weakly non-dominated faces.   
This inconsistency motivates us to develop new, consistent definitions of supported and weakly supported non-dominated points.  
 
\begin{definition}[Supported/Weakly Supported]\label{supp}
	An efficient solution is called a \emph{weakly supported efficient solution}
	if it is an optimal solution of a weighted-sum scalarization ${P}_{\lambda}$ for some weight $\lambda\in\Lambda_p^0$. Moreover, if the weight is strictly positive $\lambda\in \Lambda_p$, it is called \emph{supported efficient solution}. The corresponding image is called \emph{weakly supported} or \emph{supported (non-dominated) point}, respectively.
\end{definition}

 Note that the terminology weakly supported efficient solutions, despite the fact that these solutions are still efficient, is motivated by geometric considerations: any weakly supported non-dominated point that is not a supported non-dominated point lies on a weakly non-dominated face, but not on a non-dominated face. In contrast, supported points necessarily lie on non-dominated faces.

Let $\mathcal{Y}_{S}$ and $\mathcal{Y}_{wS}$ be the set of all supported non-dominated points and the set of all weakly supported non-dominated points, respectively. Accordingly, $\mathcal{X}_S$ ($\mathcal{X}_{wS})$ is the set of all (weakly) supported efficient solutions. We denote the cardinality of $\mathcal{X}_S$ by $|\mathcal{X}_S| = \mathcal{S}$.  

With these definitions, the following characterizations for MOCO problems hold: While weakly supported non-dominated points lie on the boundary of the upper image, supported non-dominated points lie only on the non-dominated frontier, i.\,e., only on maximally non-dominated faces. The unsupported non-dominated points lie in the interior of the upper image.

\begin{thm}
	$\Y_{wS}=\Ybd$
\end{thm}
\begin{proof}
	We follow the proofs of Theorem 3.4 and Theorem 3.5 in~\cite{ehrgott2005multicriteria} with slight modifications. To prove the equality of the two sets, we show that the subset relation holds in both directions.
	\begin{itemize}
		\item $\Y_{wS}\subseteq \Ybd$:   Let $\hat{y}\in \Y_{wS}$, i.e., there exists a $\lambda\in \Lambda^0_p$ such that $\hat{y}$ is an image of an optimal solution of~\eqref{eq:wsm} and $\hat{y}$ is non-dominated. 
		
		Suppose that $\hat{y} \notin \Ybd$, i.e., $\hat{y}\notin \partial\conv(\Y_N+\R_{\geqq}^p)$. Then $\hat{y}$ must lie in the interior of $\conv(\Y_N+\R_{\geqq}^p)$. Thus, there exists an $\varepsilon$-neighborhood $B(\hat{y}, \varepsilon)$ of $\hat{y}$, defined as $B(\hat{y}, \varepsilon) := \hat{y} + B(0, \varepsilon) \subset \conv(\Y_N+\R_{\geqq}^p)$, where $B(0, \varepsilon)$ is an open ball with radius $\varepsilon$ centered at the origin. Let $d \in \mathbb{R}^p_{>}$. 
		Then we can choose some $\alpha \in \mathbb{R}$, $0 < \alpha < \varepsilon$ such that $\alpha\, d \in B(0, \varepsilon)$. Now, $y' = \hat{y} - \alpha\, d \in \conv (\Y_N+\R_{\geqq}^p)$ with $y_k' < \hat{y}_k$ for $k = 1, \dots, p$ and
		\[
		\sum_{k=1}^p \lambda_k\, y_k' < \sum_{k=1}^p \lambda_k\, \hat{y}_k,
		\]
		because at least one of the weights $\lambda_k$ must be positive. This contradiction implies the result.

		\item $\Ybd \subseteq \Y_{wS}:$ Let $\hat{y}\in \Ybd \implies \hat{y} \in \partial\conv(\Y_N+\R_{\geqq}^p)$ and $\hat{y}$ is non-dominated, hence it exists a supporting hyperplane $\{y\in \R^p : \lambda^{\top} y = \lambda^{\top}\hat{y} \}$  with $\lambda\in \R^p\setminus\{0\}$ and it holds 
		$\lambda^{\top} y \geq \lambda^\top \hat{y}$ for all $y\in \conv(\Y_N+\R_{\geqq}^p)$.  Furthermore, by applying the separation theorem to the disjoint sets \(\{y+d\colon d\in\R^p_>\}\) and \(\{\hat{y}-d'\colon d'\in\R^p_>\}\), we obtain:
		$$
		\lambda^{\top}(y+d-\hat{y}) \geq 0 \geq \lambda^{\top}(-d') 
		$$
		for all $y\in \conv(\Y_N+\R_{\geqq}^p)$, and $d, d'\in \R^p_{>}.$ Choose $d'= e_k + \varepsilon\, e$, where $e_k$ is the $k$-th unit vector and $e=(1,\ldots, 1)^{\top}\in \R^p$. With $\varepsilon>0$ arbitrary small we see that $\lambda_k\geq 0$ for all $k= 1,\ldots, p$. Thus, $\hat{y}$ is the image of an optimal solution to~\ref{eq:wsm} with $\lambda\in \Lambda^p_0$ and per definition non-dominated. It follows $\hat{y}\in \Y_{wS}$.
	\end{itemize}
\end{proof}

Note that although the following results are stated for~\ref{eq:MOCO}, they can be directly transferred to~\ref{eq:MOILP}.

\begin{thm}\label{thm:MOCO}
    For \ref{eq:MOCO} problems, it holds $\Y_S = \Y_{SNF}.$
\end{thm}

\begin{proof}
	$\Y_{S} \subseteq \Y_{SNF}:$ Let $\hat{y}\in \Y_{S}=\Y_{S\lambda}$, i.e., $\hat{y}$ is the image of an optimal solution of $P_\lambda$ with $\lambda \in \Lambda_p$. Since $\conv (\Y_N)$ is a nonempty convex set, $\hat{y}$ is a minimal element of $\conv(\Y_N)$ induced of the cone $-\R^p_{\geqq}$ and hence it holds $((\hat{y}-\R^p_{\geqq})\backslash{\hat{y}})\cap \conv(\Y_N) = \varnothing$ \citep{boyd2004convex}. Consequently, it follows that $\hat{y}\in \Y_{SNF}$.

	$\Y_{S} \supseteq \Y_{SNF}:$  Any point $\hat{y}\in \Y_{SNF}$ is a minimal element of the set $\conv(\Y_N)$ induced by the cone $-\R_{\geqq}^p$. Thus $\hat{y}$ is a non-dominated point for the set $\conv(\Y_N)$. If we define the relaxation MOLP of the given~\ref{eq:MOCO} problem by
	\[
	\min_{x\in \conv(\X)} C\, x, 
	\]
	it follows that $\hat{y}$, is an image of an optimal solution for the reduced MOLP. According to~\Cref{thm:iser}, there is a  $\lambda\in \Lambda_p$ such that $\hat{y}$ is the image of an optimal solution to~\ref{eq:wsm} and hence $\hat{y}\in \Y_S$. 
\end{proof}

Consequently, Definition~\ref{def:supp}~(1) includes the weakly supported non-dominated points for general \ref{eq:MOP} problems, while Definitions~(3) and (4) only contain the supported non-dominated points (in the discrete case). Hence, in the discrete case, the weakly supported non-dominated points lie on the boundary of the upper image while the supported ones lie only on the non-dominated frontier, 
i.\,e., on maximally non-dominated faces. 

\begin{figure}[h]
\centering
\begin{minipage}[t]{0.45\textwidth}
    \centering
    \begin{tikzpicture}[]
        \tdplotsetmaincoords{120}{60}
        \tikzstyle{vertex}=[circle,fill=black,draw=black,minimum size=3pt,inner sep=0]
        \tikzstyle{vertex2}=[circle,fill=red,draw=red,minimum size=3pt,inner sep = 0]
        \tikzstyle{vertex4}=[diamond,fill=black,draw=black,minimum size=4pt,inner sep = 0]
        \newcommand{\Cross}{$\mathbin{\tikz [x=1ex,y=1ex,line width=.1ex, black] \draw (0,0) -- (1,1) (0,1) -- (1,0);}$}%
        \newcommand{\Crossa}{$\mathbin{\tikz [x=1ex,y=1ex,line width=.1ex, blue] \draw (0,0) -- (1,1) (0,1) -- (1,0);}$}%
        \tikzstyle{N_point}=[draw, cross out,scale=.5,thick]
        \tikzstyle{U_point}=[fill, circle,scale=.5]
        \tikzstyle{non_maximal}=[draw,circle,scale=.5, thick]
        \tikzstyle{U_point_act}=[fill, red, circle,scale=.5]
        \tikzstyle{U_point_visited}=[draw, strike out,scale=1,thick]
        \tikzstyle{S_zone}=[fill, blue!30]
        \tikzstyle{S_zone_rect}=[S_zone, rectangle, scale=.9]
        \tikzstyle{D_zone}=[fill, mathdarkblue!30]
        \tikzstyle{D_zone_rect}=[D_zone, rectangle, scale=.9]
        \tikzstyle{lface}=[mathdarkblue!20,semitransparent]
        \tikzstyle{sface}=[mathdarkblue!40,semitransparent]
        \tikzstyle{sfaceu}=[unigreen!80, semitransparent]
        \tikzstyle{dface}=[mathdarkblue!60,semitransparent]
        \tikzstyle{VSplit}=[
            thick,
            draw=mathdarkblue!50!black!50,
            top color=mathdarkblue!10,
            bottom color=mathdarkblue!80!black!20,
            minimum size=1em,
            rectangle split,
            %rectangle split horizontal=false,
            rectangle split ignore empty parts] % requires library shapes.multipart

        \providecommand*\ColVec[3]{
            (#1,#2,#3)
        }

        \providecommand*{\VecNode}[4]{node[anchor=west] {$#1=(#2,#3,#4)$}}

        \begin{scope}[tdplot_main_coords, scale=.35]
            % Feasible points
            \foreach \x/\y/\z/\name in {3/8/8/z1, 4/6/5/z2}
                \coordinate (\name) at (\x,\y,\z);

            % Local upper bounds
            \foreach \x/\y/\z/\name in {10/10/10/M,3/10/10/u1, 10/10/8/u3, 10/8/10/u2, 
            6/10/7/u31, 6/5/10/u21, 10/2/10/u22, 10/5/4/u23, 10/10/4/u33, 10/2/7/u32}
                \coordinate (\name) at (\x,\y,\z);

            {
                % \node[U_point, label={[]below:$0$}] at (0,0,0) {};
            }

            \draw[] (0,0,0) -- (10,0,0);
            \draw[] (0,0,0) -- (0,10,0);
            \draw[] (0,0,0) -- (0,0,10); 
            \draw[->] (10,0,0) node[anchor=north]{} -- (13,0,0) node[anchor=south]{{$y_3$}};
            \draw[->] (0,10,0) -- (0,13,0) node[anchor=north west]{{$y_1$}};
            \draw[->] (0,0,10) -- (0,0,13) node[anchor=south]{{$y_2$}};

            \fill[DodgerBlue4,semitransparent] (5,0,0) -- (5,0,12) -- (5,12,12) -- (5,12,0);
            \node[] at (5,9.8,10.5){$y_3 =5$};
            \node[DodgerBlue4!50!black,vertex, label={[label distance=-5pt]above right:{}}] at (5,11,1) {};
            \node[DodgerBlue4!50!black,vertex, label={[label distance=-5pt]above right:{}}] at (5,1,10) {};
            \node[DodgerBlue4!50!black,vertex, label={[label distance=-5pt]above right:{}}] at (5,3,4) {};
            \fill[DeepPink4,semitransparent] (1,0,0) -- (1,0,12) -- (1,12,12) -- (1,12,0);

            \node[] at (1,10,10 ){$y_3 =1$};

            \node[vertex, label={[label distance=-5pt]above left:{$y^1$}}] at (1,2,9) {};
            \node[vertex, label={[label distance=-5pt]above left:{$y^2$}}] at (1,3,6) {};
            \node[N_point, label={[label distance=-5pt]above right:{$y^4$}}] at (1,6,5) {};
            \node[vertex, label={[label distance=-5pt]above right:{$y^3$}}] at (1,8,3) {};
        \end{scope}    
    \end{tikzpicture}
\end{minipage}%
\hspace{0.05\textwidth}
\begin{minipage}[t]{0.45\textwidth}
    \centering
    \begin{tikzpicture}[]
        \tdplotsetmaincoords{0}{0}
        \tikzstyle{vertex}=[circle,fill=black,draw=black,minimum size=3pt,inner sep=0]
        \tikzstyle{vertex2}=[circle,fill=red,draw=red,minimum size=3pt,inner sep = 0]
        \tikzstyle{N_point}=[draw, cross out,scale=.5,thick]
        \tikzstyle{U_point}=[fill, circle,scale=.5]
        \tikzstyle{non_maximal}=[draw,circle,scale=.5, thick]
        \tikzstyle{U_point_act}=[fill, red, circle,scale=.5]
        \tikzstyle{U_point_visited}=[draw, strike out,scale=1,thick]
        \tikzstyle{S_zone}=[fill, blue!30]
        \tikzstyle{S_zone_rect}=[S_zone, rectangle, scale=.9]
        \tikzstyle{D_zone}=[fill, mathdarkblue!30]
        \tikzstyle{D_zone_rect}=[D_zone, rectangle, scale=.9]
        \tikzstyle{lface}=[unigreen!20, semitransparent]
        \tikzstyle{sface}=[unigreen!80, semitransparent]
        \tikzstyle{sfaced} = [mathdarkblue!40,semitransparent]
        \tikzstyle{dface}=[black!60, semitransparent]
        \tikzstyle{VSplit}=[
            thick,
            draw=mathdarkblue!50!black!50,
            top color=mathdarkblue!10,
            bottom color=mathdarkblue!80!black!20,
            minimum size=1em,
            rectangle split,
            rectangle split ignore empty parts] % requires library shapes.multipart
        \providecommand*\ColVec[3]{
            (#1,#2,#3)
        }
        \providecommand*{\VecNode}[4]{node[anchor=west] {$#1=(#2,#3,#4)$}}

        \begin{scope}[tdplot_main_coords, scale=.35]
          
	 {
		% \node[U_point, label={[]below:$0$}] at (0,0) {};
	}

	 {
	
        %\node[N_point, label={[label distance=-5pt]above left:{$z^1$}}] at (3,8,8) {};
	}

	{
	%	\node[U_point, label={[label distance=-5pt]below right:$u^{2}$}] at (u2) {};
	%	\node[U_point, label={[label distance=-5pt]below right:$u^{3}$}] at (u3) {};
	}
	{
		%\node[U_point, label={[label distance=-5pt]below right:$u^{33}$}] at (u33) {};
		%\node[U_point, label={[label distance=-5pt]below right:$u^{31}$}] at (u31) {};
	}

	{
%		\node[U_point, label={below:$u^{23}$}] at (u23) {};
	}
	{
%		\node[U_point, label={[label distance=-5pt]below right:$u^{32}$}] at (u32) {};
	}

	{
%		\node[U_point_visited] at (u23) {};
%		\node[U_point_visited] at (u32) {};
	}

	%\visible<7> {
		%% Edges associated to u3
		%\draw[thick] (u21) -- (3,5,10);
		%\draw[thick] (u21) -- (6,2,10);
		%\draw[thick] (u21) -- (6,5,7);
		%%\draw[thick, dashed] (6,5,7);
	%}

        \fill[DeepPink4,semitransparent] (0,0) -- (0,12) -- (12,12) -- (12,0);
        \draw[] (0,0) -- (10,0);
        \draw[] (0,0) -- (0,10);
	\draw[->] (10,0) node[anchor=north]{} -- (13,0) node[anchor=north]{{$y_1$}};
	\draw[->] (0,10) -- (0,13) node[anchor=north east]{{$y_2$}};

        \draw[thick] (2,9) -- (3,6) -- (8,3);
             \draw[dotted,thick] (2,9) -- (2,12);
             \draw[dotted,thick] (8,3) -- (12,3);
        %\fill[lightgray, semitransparent] (2,9) -- (2,12) -- (12,12) -- (12,3) -- (8,3) -- (3,6) -- (2,9);

           % \fill[blue,semitransparent] (5,0,0) -- (5,0,12) -- (5,12,12) -- (5,12,0);

        %\node[blue,N_point, label={[label distance=-5pt]above right:{}}] at (5,11,1) {};
          %\node[blue,N_point, label={[label distance=-5pt]above right:{}}] at (5,6,6) {};
         %  \node[blue,N_point, label={[label distance=-5pt]above right:{}}] at (5,1,10) {};
        %\node[blue,N_point, label={[label distance=-5pt]above right:{}}] at (5,3,4) {};

        \node[] at ( 10 ,11 ){${y_3 =1} $};
        % \node[] at (5,10,10 ){${\color{blue}c_3 =5} $};

        \node[vertex, label={[label distance=-5pt]above right:{$y^1$}}] at (2,9) {};
        \node[vertex, label={[label distance=-5pt]above right:{$y^2$}}] at (3,6) {};
        \node[N_point, label={[label distance=-5pt]above right:{$y^4$}}] at (6,5) {};
        \node[vertex, label={[label distance=-5pt]above right:{$y^3$}}] at (8,3) {};
        %\draw (2,9) -- (3,6) -- (8,3);
        \end{scope}
    \end{tikzpicture}
\end{minipage}
 \caption{Outcome space with non-dominated points $y^1,\ldots,y^4$, each with the value of 1 in the third component and the two-dimensional projection of the plane $c_3=1$.}\label{fig:2dvs3d}
\end{figure}
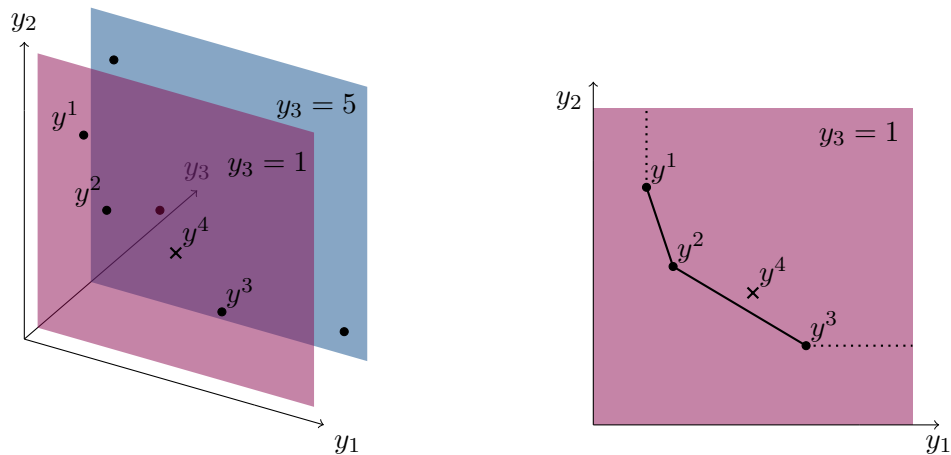

\begin{example}\label{ex:weaklysupp}
To illustrate the geometrical properties of weakly supported non-dominated points, consider the outcome vectors in~\Cref{fig:2dvs3d}. The outcome vectors are partitioned into two layers based on their value in the third component. The non-dominated points \( y^1 = (2, 9, 1), y^2 = (3, 6, 1), y^3 = (8, 3, 1), \) and \( y^4 = (6, 5, 1) \) all share a minimum value of \( c_3 = 1 \) in their third component (pink layer), while there are other non-dominated points having a value of 5 in the third component (blue layer). 

Among them, only \( y^1, y^2, \) and \( y^3 \) lie on the non-dominated frontier and can be obtained as optimal solutions of a weighted-sum scalarization with \( \lambda \in \Lambda_d \). In contrast, \( y^4 \) does not lie on the non-dominated frontier but is still on the boundary of the upper image (on a weakly non-dominated face). It can be obtained as an optimal solution of a weighted-sum scalarization where some weights are zero, e.g., obtained as an optimal solution of $P_\lambda$ with $\lambda^{\top}=(0,0,1)^{\top}$, making \( y^4 \) a weakly supported but not a supported non-dominated vector.
\end{example}

 \begin{thm}\label{thm:supp_not_equivalent}
 Let $\mathcal{Y}_{S}$ and $\mathcal{Y}_{wS}$ be the sets of supported and weakly supported non-dominated points of a MOCO problem, respectively. Then, $\mathcal{Y}_{S} \subseteq \mathcal{Y}_{wS}$ and there exist instances where $\mathcal{Y}_{S} \subset \mathcal{Y}_{wS}$, i.\,e.,  the set of supported non-dominated points is a proper subset of the set of the weakly supported non-dominated points.  \end{thm}

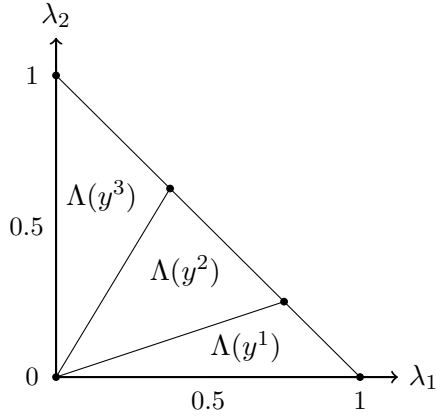
\begin{figure}[h]
\centering
\begin{tikzpicture}
\tikzstyle{vertex}=[circle,fill=black,draw=black,minimum size=2.5pt,inner sep=0]
%uncomment if require: \path (0,353); %set diagram left start at 0, and has height of 353

%Shape: Axis 2D 
\draw [<->,thick] (0,4.5) node (yaxis) [above] {$\lambda_2$}
        |- (4.5,0) node (xaxis) [right] {$\lambda_1$};
        
\draw (0,4) node (a_1) [xshift=-9pt] {\small$1$} -- (4,0) node (a_2) [yshift=-9pt] {\small$1$};
\draw (0,2) node (a_2) [xshift=-11pt] {\small$0.5$} -- (0,4) coordinate (a_3);
\draw (0,0) node (a_4) [xshift=-9pt] {\small$0$} -- (3,1) node (a_5) {};
\draw (2,0) node (a_6) [yshift=-9pt] {\small$0.5$} -- (4,0) coordinate (a_7);
\draw (0,0)--(3/2,5/2);

\node[vertex] at (0,0) (int1) {};
\node[vertex] at (4,0) (int1) {};
\node[vertex] at (0,4) (int1) {};
\node[vertex] at (3/2,5/2) (int1) {};
\node[vertex] at (3,1) (int1) {};

\node at (0.6,2.4) (int1) {$\Lambda(y^3)$};
\node at (1.7,1.4) (int1) {$\Lambda(y^2)$};
\node at (2.5,0.4) (int1) {$\Lambda(y^1)$};
\end{tikzpicture}
 \caption{Projected weight space decomposition to the upper image of $\Y = \{(2, 9, 1)^\top, (3, 6, 1)^\top, (8, 3, 1)^\top,(6, 5, 1)^\top\}$ with $\lambda_3\coloneqq 1-\lambda_2-\lambda_1.$ The set of weighting vectors associated with a point $y\in \mathcal{Y}$ is given by 
$\Lambda(y)\coloneqq \{ \lambda \in \Lambda^0_{p}\colon \lambda^\top  y \leq \lambda^\top  y^{\prime} \text{ for all } y^{\prime} \in \mathcal{Y^{\geq}}\}.$ For a comprehensive overview of the weight space decomposition, we refer to~\cite{Przy10}.}
 \label{fig:wsd}
\end{figure}

\begin{proof}
	The inclusion $\mathcal{Y}_{S} \subseteq \mathcal{Y}_{wS}$ holds per definition.
	Consider the outcome set $\Y = \{y^1,y^2,y^3,y^4\}$ with  \( y^1 = (2, 9, 1), y^2 = (3, 6, 1), y^3 = (8, 3, 1), \) and \( y^4 = (6, 5, 1) \), which represents the pink layer of the outcome set in~\Cref{ex:weaklysupp}, displayed in $\Cref{fig:2dvs3d}$.  It is straightforward to construct an artificial MOILP problem with this outcome set. 
	
	All points $y^i$ with $i\in \{1,\ldots, 4\}$  are weakly supported non-dominated points since their preimages are optimal solutions of a weighted sum problem $P_\lambda$ with $\lambda= (0 , 0, 1)^\top $. 
	Since $y^1$,$y^2$ and $y^3$ are also supported non-dominated points, since their preimages are optimal solutions of the respective weighted sum problems $P_{\lambda^1}$ with $\lambda^1= (0.7 , 0.1, 0.2)^\top $, $P_{\lambda^2}$ with $\lambda^2= (0.4 , 0.4, 0.2)^\top $, and $P_{\lambda^3}$ with $\lambda^3= (0.1 , 0.8, 0.1)^\top $, respectively. The corresponding \emph{(projected) weight space decomposition} is illustrated in~\Cref{fig:wsd}.  Note that, for example, $(6,4.2,1)^{\top} \in \conv(\mathcal{Y}_N) \cap (y^4 -\R^p_{\geqq})$ implies $y^4 \notin \{y \in \conv(\mathcal{Y}_N) \colon \conv(\mathcal{Y}_N) \cap (y -\R^p_{\geqq}) = \{y\}\}$. By \Cref{thm:MOCO}, it follows that $y_4 \notin \Y_{S}$, although $y^4\in \Y_{wS}$.  
	
	Thus, the set of weakly supported non-dominated points is $\mathcal{Y}_{wS}=\{y^1,y^2,y^3,y^4\}$, while the set of supported non-dominated points is $\mathcal{Y}_{S}=\{y^1,y^2,y^3 \}$. Hence, in this example, the set of supported non-dominated points is a proper subset of the set of weakly supported non-dominated points $\mathcal{Y}_{S} \subset \mathcal{Y}_{wS}$.
\end{proof}

Note that in the bi-objective case $\mathcal{Y}_{S} = \mathcal{Y}_{wS}$ as the following lemma shows.
\begin{lem}
  Every weakly supported non-dominated point of a biobjective integer optimization problem is supported non-dominated.
\end{lem}
\begin{proof}
  Let \(\bar{y}=f(\bar{x})\in\R^2\) be a weakly supported but not supported non-dominated point of a bi-objective integer optimization problem and \(\bar{x}\) the corresponding preimage. Then there is a weighting vector \(\lambda\in\Lambda_0\) such that \(\bar{x}\) is an optimal solution of the weighted sum problem \(P_{\lambda}\). Since \(\bar{y}\) is not supported non-dominated, one of the components of \(\lambda\) must be zero, w.l.o.g.\ let \(\lambda=(1,0)^\top\).

  Since \(y\) is non-dominated, there does not exist a feasible outcome vector \(y=(y_1,y_2)^\top\), with \(y_1=\bar{y}_1\) and \(y_2<\bar{y}_2\). Thus, \(\bar{x}\) is also an optimal solution of the weighted sum problem \(P_{\lambda'}\) with \(\lambda'=(1-\varepsilon,\varepsilon)^\top \) for $\varepsilon>0$ sufficiently small, which makes \(\bar{y}\) supported non-dominated.
\end{proof}

The clear distinction between the sets of supported non-dominated and weakly supported non-dominated solutions  (\Cref{supp}) is also necessary as the corresponding
problems may differ in their output time complexity.  For instance, in the case of the minimum cost flow problem, it can be shown that supported efficient solutions can be determined in \emph{output-polynomial time}, while this is not the case for \emph{weakly supported solutions}, unless $\mathbf{P}= \mathbf{NP}$~\citep{konen2023outputpolynomial}.

\begin{thm}
    The determination of all weakly-supported solutions for a \ref{eq:MOCO} problem with $p+1$ objectives is as hard as the determination of all non-dominated points for a \ref{eq:MOCO} with $p$ objectives. 
\end{thm}
\begin{proof}
    Assume an algorithm exists to determine all weakly
    supported non-dominated points for a given \ref{eq:MOCO} problem with $p+1$ objectives. Let $M_p$ be a \ref{eq:MOCO} with $p$~objectives. Suppose we add an artificial objective $c^{p+1} =0$ to our \ref{eq:MOCO} problem and denote it by $M_{p+1}$. In that case, we obtain a weakly efficient facet for $M_{p+1}$, where all non-dominated points for $M_p$  are \weakly supported non-dominated points for $M_{p+1}$. Therefore, even the unsupported non-dominated points for $M_p$ are \weakly supported non-dominated points for $M_{p+1}$ since they are part of the boundary of the upper image for $M_{p+1}$. Consequently, any algorithm that can determine all weakly supported non-dominated points for $M_{p+1}$  can determine all non-dominated points for $M_p$.
\end{proof}

While the above characterization of supported points holds for \ref{eq:MOCO} problems,it does not extend to general \ref{eq:MOP} problems, and it can be shown that in the case of general \ref{eq:MOP} problems it may hold $\Y_{S} \subset \Y_{SNF}$.  This distinction arises because the weighted sum scalarization for \(\lambda\in\Lambda_p\) method is only capable of identifying \emph{properly supported non-dominated points} in the sense of Geoffrion~\citep{GEOFFRION1968618}, as discussed in~\cite{ehrgott2005multicriteria}. Note that for \ref{eq:MOCO} problems every efficient solution is also properly efficient \citep{ehrgott2005multicriteria}. However, for general \ref{eq:MOP} problems, the set of properly supported non-dominated points may form a strict subset of the non-dominated points on the non-dominated frontier, highlighting the limitations of the weighted sum approach in capturing the entire set of supported solutions in general MOO problems~\citep{ehrgott2005multicriteria}.  The following example, similar to the ones in \cite{boyd2004convex} or \cite{ehrgott2005multicriteria}, illustrates this. 

\begin{figure}
    \centering
    \begin{tikzpicture}[scale=1.5]
        % Axes
        \draw[<->] (0,3) node[left] {$f_2$} -- (0,0) -- (3,0) node[below] {$f_1$};
        % \draw[->]  -- (0,3) ;
        
        % Shading for y^1 + R^2_{>=}
        \fill[black!20,opacity=0.35] (1,2) -- (-0.3,2) -- (-0.3,-0.3) -- (1,-0.3) -- cycle ;
        \draw[darkgray]  (1,-0.33) -- (1,2) -- (-0.33,2) ;
        
        % Shading for y^2 + R^2_{>=}
        \fill[black!20,opacity=0.35] (2,1) -- (2,-0.3) -- (-0.3,-0.3) -- (-0.3,1) -- cycle ;
        \draw[darkgray]  (-0.33,1) -- (2,1) -- (2,-0.33) ;
        
        % Labels for y^1 + R^2_{>=} and y^2 + R^2_{>=}
        \node[black] at (1.5,0.5) {$y^2 - \mathds{R}^2_{\geqq}$};
        \node[black] at (0.5,1.5) {$y^1 - \mathds{R}^2_{\geqq}$};
        \node[DodgerBlue4] at (2,2) {$\Y$};
        % Feasible region
        \fill[DodgerBlue4,opacity=0.4] (2,2) circle (1cm);
        
        % Points y^1 and y^2
        \fill[] (1,2) circle (0.05) node[above] {$y^1=(\begin{smallmatrix}1\\2\end{smallmatrix})$};
        % \draw[red] (1,2) circle (0.1);
        \fill[] (2,1) circle (0.05) node[below right] {$y^2=(\begin{smallmatrix}2\\1\end{smallmatrix})$};
        % \draw[red] (2,1) circle (0.1);
        
        % Lower left boundary (thick curve)
        \draw[thick] (1,2) arc[start angle=180,end angle=270,radius=1cm];
        %\draw[thick] (2,1) arc[start angle=0,end angle=45,radius=2cm];
    \end{tikzpicture}
     \caption{Illustration of~\Cref{ex:MOO} highlighting the points $y^1,y^2$ the points  that contains the points that are both in $\Y_{SNF}$ and $\Y_{S\lambda}$ (bold curve), which together with $y^1$ and $y^2$ gives the non-dominated frontier.  }
     \label{fig:MOO_gen}
 \end{figure}
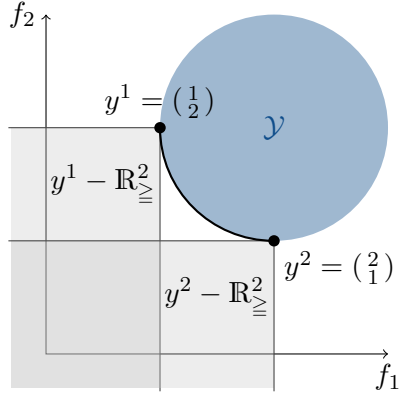

\begin{example}\label{ex:MOO}
    Consider the following \ref{eq:MOP} problem 
    \begin{align*}
        \min \quad& x_1+2 \\ 
        \min \quad& x_2 +2\\
        \text{s.t.} \quad & x_1^2 + x_2^2 \leq 1 \\
    \end{align*}
    Consider the points $y^1= (1,2)$ and $y^2=(2,1)$, which can be obtained with $\lambda^1=(0,1)$ and $\lambda^2=(1,0)$, respectively. It holds $y^1,y^2\in \Y_{SNF}$ but $y^1$ and $y^2$ are not optimal solutions for $P_\lambda$ for any $\lambda\in \Lambda$~\citep{boyd2004convex,ehrgott2005multicriteria}. Thus, $y^1,y^2\notin \Y_{S}$. However, any point on the left lower boundary without the points $(1,2)$ and $(2,1)$ are both contained in $\Y_{SNF}$ and $\Y_{S}$. Therefore, in this example,  $\Y_{S}\subset \Ysnf$. The example is illustrated in~\Cref{fig:MOO_gen}
\end{example}

\begin{thm}
    For general~\ref{eq:MOP} problems it holds $\Y_{S\lambda}\subseteq \Y_{SNF}$ and  there exist instances where  $\Y_{S\lambda}\subset \Y_{SNF}$.
\end{thm}
\begin{proof}
    $\Ywsm \subseteq \Ysnf$ follows directly from the proof of~\Cref{thm:MOCO}.
    Furthermore, as demonstrated in~\Cref{ex:MOO}, there exist instances where $\Y_{S\lambda} \subset \Y_{SNF}$.
\end{proof}

The distinction between supported and weakly supported is sufficient for \ref{eq:MOCO} and \ref{eq:MOILP} problems. However, in the context of general \ref{eq:MOP}, it may be worthwhile to introduce a finer classification, distinguishing between properly supported, supported, and weakly supported non-dominated points. While this work primarily focuses on the supportedness in \ref{eq:MOCO} problems, a comprehensive analysis of supportedness in general \ref{eq:MOP} problems remains an important area for future research. Notably, recent advancements in this direction have been discussed in \cite{ChlumskyHarttmann2025}.

\section{Conclusion}\label{chapt:concl}

 Some previous papers use inconsistent characterizations of supported non-dominated points for \ref{eq:MOCO} problems. Through counterexamples and theoretical analysis, this paper proves that these definitions, while being equivalent in the context of~\ref{eq:MOLP}, diverge in~\ref{eq:MOCO} problems, yielding distinct sets of supported non-dominated points with differing structural and computational properties. This emphasizes the need for precise and consistent definitions.  Hence, this paper proposes the definition of weakly supported non-dominated points and establishes a clear distinction between weakly supported and supported non-dominated points. This is particularly important for \ref{eq:MOCO} problems with more than two objectives.

 Using these refined definitions, the following characterizations for MOCO problems hold: While the weakly supported non-dominated points lie on the boundary of the upper image, the supported non-dominated points lie only on the non-dominated frontier, i.\,e., only on  maximally non-dominated faces, and the unsupported non-dominated points lie in the interior of the upper image.  However, these characterizations do not extend directly to general MOO problems. Hence, future research should aim to extend the characterization of supportedness beyond MOCO problems to general MOO contexts, addressing the gaps highlighted by this study. An approach towards a categorization of supportedness in the context of general multi-objective optimization has been discussed in~\cite{ChlumskyHarttmann2025}.

Unsupported solutions may be reasonable compromise solutions and should thus not be neglected completely. Note that the difficulty in computing unsupported solutions arises in many integer and combinatorial optimization problems and is one reason for their computational complexity, in general~\citep{ehrgott00hard, figueira17easy}. One way to overcome this computational burden---at least to a certain degree---could be to determine unsupported solutions only in regions of the Pareto front that are not well represented by the set of supported non-dominated points.

\bigskip

\section*{Acknowledgements.}  
The authors thankfully acknowledge the financial support of the Deutsche Forschungs\-gemeinschaft (DFG, German Research Foundation), project number 441310140.

%\begin{thebibliography}{00}
\providecommand{\natexlab}[1]{#1}
	\providecommand{\url}[1]{{#1}}
	\providecommand{\urlprefix}{URL }
	\expandafter\ifx\csname urlstyle\endcsname\relax
	\providecommand{\doi}[1]{DOI~\discretionary{}{}{}#1}\else
	\providecommand{\doi}{DOI~\discretionary{}{}{}\begingroup
		\urlstyle{rm}\Url}\fi
	\providecommand{\eprint}[2][]{\url{#2}}

\bibliographystyle{apalike}

\begin{thebibliography}{}
	
	\bibitem[Argyris et~al., 2011]{argyris2011identifying}
	Argyris, N., Figueira, J.~R., and Morton, A. (2011).
	\newblock Identifying preferred solutions to multi-objective binary
	optimisation problems, with an application to the multi-objective knapsack
	problem.
	\newblock {\em Journal of Global Optimization}, 49:213-235.
	
	\bibitem[Benson, 1998]{Benson1998AnOA}
	Benson, H.~P. (1998).
	\newblock An outer approximation algorithm for generating all efficient extreme
	points in the outcome set of a multiple objective linear programming problem.
	\newblock {\em Journal of Global Optimization}, 13:1-24.
	
	\bibitem[Boyd and Vandenberghe, 2004]{boyd2004convex}
	Boyd, S. and Vandenberghe, L. (2004).
	\newblock {\em Convex optimization}.
	\newblock Cambridge university press.
	
	\bibitem[B\"okler, 2018]{boeklerthesis}
	B\"okler, F. (2018).
	\newblock {\em Output-sensitive Complexity of Multiobjective Combinatorial
		Optimization Problems with an Application to the Multiobjective Shortest Path
		Problem}.
	\newblock PhD thesis.
	
	\bibitem[B\"okler et~al., 2024]{Bokler2024}
	B\"okler, F., Parragh, S.~N., Sinnl, M., and Tricoire, F. (2024).
	\newblock An outer approximation algorithm for generating the
	edgeworth-pareto hull of multi-objective mixed-integer linear programming
	problems.
	\newblock {\em Mathematical Methods of Operations Research}, 100(1):263-290.
	
	
	\bibitem[Chlumsky-Harttmann, 2025]{ChlumskyHarttmann2025}
	Chlumsky-Harttmann, F. (2025).
	\newblock {\em Robust Multi-Objective Optimization: Analysis and Algorithmic
		Approaches}.
	\newblock doctoralthesis, Rheinland-Pfälzische Technische Universität
	Kaiserslautern-Landau.
	
	\bibitem[Correia et~al., 2021]{Correia2021}
	Correia, P., Paquete, L., and Figueira, J.~R. (2021).
	\newblock Finding multi-objective supported efficient spanning trees.
	\newblock {\em Computational Optimization and Applications}, 78(2):491-528.
	
	\bibitem[da~Silva and Cl\'i­maco, 2007]{silva07note}
	da~Silva, C.~G. and Cl\'i­maco, J. (2007).
	\newblock A note on the computation of supported non-dominated solutions in the
	bi-criteria minimum spanning tree problem.
	\newblock {\em Journal of Telecommunications and Information Technology}, page
	11-15.
	
	\bibitem[Dai and Charkhgard, 2018]{dai2018two}
	Dai, R. and Charkhgard, H. (2018).
	\newblock A two-stage approach for bi-objective integer linear programming.
	\newblock {\em Operations Research Letters}, 46(1):81-87.
	
	\bibitem[{Edwin Romeijn} and Smith, 1999]{EDWIN99}
	{Edwin Romeijn}, H. and Smith, R.~L. (1999).
	\newblock Parallel algorithms for solving aggregated shortest-path problems.
	\newblock {\em Computers \& Operations Research}, 26(10):941-953.
	
	\bibitem[Ehrgott, 2000]{ehrgott00hard}
	Ehrgott, M. (2000).
	\newblock Hard to say it's easy - four reasons why combinatorial
	multiobjective programmes are hard.
	\newblock In Haimes, Y.~Y. and Steuer, R.~E., editors, {\em Research and
		Practice in Multiple Criteria Decision Making}, volume 487 of {\em Lecture
		Notes in Economics and Mathematical Systems}, page 69-80. Springer Berlin
	Heidelberg.
	
	\bibitem[Ehrgott, 2005]{ehrgott2005multicriteria}
	Ehrgott, M. (2005).
	\newblock {\em Multicriteria optimization}, volume 491.
	\newblock Springer Science \& Business Media.
	
	\bibitem[Eus\'ebio and Figueira, 2009a]{EUSEBIO20092554}
	Eus\'ebio, A. and Figueira, J.~R. (2009a).
	\newblock Finding non-dominated solutions in bi-objective integer network flow
	problems.
	\newblock {\em Computers \& Operations Research}, 36(9):2554-2564.
	
	\bibitem[Eus\'ebio and Figueira, 2009b]{Eusebio09}
	Eus\'ebio, A. and Figueira, J.~R. (2009b).
	\newblock Finding non-dominated solutions in bi-objective integer network flow
	problems.
	\newblock {\em Computers \& Operations Research}, 36:2554-2564.
	
	\bibitem[Eus\'ebio and Figueira, 2009c]{EUSEBIO200968}
	Eus\'ebio, A. and Figueira, J.~R. (2009c).
	\newblock On the computation of all supported efficient solutions in
	multi-objective integer network flow problems.
	\newblock {\em European Journal of Operational Research}, 199(1):68-76.
	
	\bibitem[Figueira et~al., 2017]{figueira17easy}
	Figueira, J.~R., Fonseca, C.~M., Halffmann, P., Klamroth, K., Paquete, L.,
	Ruzika, S., Schulze, B., Stiglmayr, M., and Willems, D. (2017).
	\newblock Easy to say they're hard, but hard to see they're easy - towards
	a categorization of tractable multiobjective combinatorial optimization
	problems.
	\newblock {\em Journal of Multi-Criteria Decision Analysis}, 24(1-2):82-98.
	
	\bibitem[Gandibleux et~al., 2001]{gandibleux01}
	Gandibleux, X., Morita, H., and Katoh, N. (2001).
	\newblock The supported solutions used as a genetic information in a population
	heuristic.
	\newblock In Zitzler, E., Thiele, L., Deb, K., {Coello Coello}, C.~A., and
	Corne, D., editors, {\em Evolutionary Multi-Criterion Optimization}, page
	429-442, Berlin, Heidelberg. Springer Berlin Heidelberg.
	
	\bibitem[Gandibleux et~al., 2003]{gandibleux2003use}
	Gandibleux, X., Morita, H., and Katoh, N. (2003).
	\newblock Use of a genetic heritage for solving the assignment problem with two
	objectives.
	\newblock In {\em International Conference on Evolutionary Multi-Criterion
		Optimization}, page 43-57. Springer.
	
	\bibitem[Geoffrion, 1968]{GEOFFRION1968618}
	Geoffrion, A.~M. (1968).
	\newblock Proper efficiency and the theory of vector maximization.
	\newblock {\em Journal of Mathematical Analysis and Applications},
	22(3):618-630.
	
	\bibitem[Hamacher et~al., 2007]{Hamacher2007}
	Hamacher, H.~W., Pedersen, C.~R., and Ruzika, S. (2007).
	\newblock Multiple objective minimum cost flow problems: A review.
	\newblock {\em European Journal of Operational Research}, 176(3):1404-1422.
	
	\bibitem[Isermann, 1974]{Iser74m}
	Isermann, H. (1974).
	\newblock Technical note - proper efficiency and the linear vector maximum
	problem.
	\newblock {\em Operations Research}, 22(1):189-191.
	
	\bibitem[Jesus, 2015]{jesus2015implicit}
	Jesus, A. D. B.~d. (2015).
	\newblock Implicit enumeration for representation systems in multi-objective
	optimization.
	\newblock Master's thesis.
		
	\bibitem[K\"onen and Stiglmayr, 2023]{konen2023outputsensitive}
	K\"onen, D. and Stiglmayr, M. (2023).
	\newblock Output-sensitive complexity of multi-objective integer network flow
	problems. \emph{arXiv preprint: 2312.01786}.
	
	\bibitem[K\"onen and Stiglmayr, 2025]{konen2023outputpolynomial}
	K\"onen, D. and Stiglmayr, M. (2025).
	\newblock An output-polynomial time algorithm to determine all supported
	efficient solutions for multi-objective integer network flow problems.
	\newblock {\em Discrete Applied Mathematics}, 376:1–14.
	
	\bibitem[Liefooghe et~al., 2014]{liefooghe2014hybrid}
	Liefooghe, A., Verel, S., and Hao, J.-K. (2014).
	\newblock A hybrid metaheuristic for multiobjective unconstrained binary
	quadratic programming.
	\newblock {\em Applied Soft Computing}, 16:10-19.
	
	\bibitem[Liefooghe et~al., 2015]{liefooghe2015experiments}
	Liefooghe, A., Verel, S., Paquete, L., and Hao, J.-K. (2015).
	\newblock Experiments on local search for bi-objective unconstrained binary
	quadratic programming.
	\newblock In {\em Evolutionary Multi-Criterion Optimization: 8th International
		Conference, EMO 2015, GuimarÃ£es, Portugal, March 29 April 1, 2015.
		Proceedings, Part I 8}, page 171-186. Springer.
	
	\bibitem[Medrano and Church, 2014]{church14}
	Medrano, F.~A. and Church, R.~L. (2014).
	\newblock Corridor location for infrastructure development: A fast bi-objective
	shortest path method for approximating the pareto frontier.
	\newblock {\em International Regional Science Review}, 37(2):129-148.
	
	\bibitem[Medrano and Church, 2015]{Church2015}
	Medrano, F.~A. and Church, R.~L. (2015).
	\newblock A parallel computing framework for finding the supported solutions to
	a biobjective network optimization problem.
	\newblock {\em Journal of Multi-Criteria Decision Analysis}, 22(5-6):244-259.
	
	\bibitem[Nemhauser and Wolsey, 1999]{nem99}
	Nemhauser, G.~L. and Wolsey, L.~A. (1999).
	\newblock {\em Integer and combinatorial optimization}.
	\newblock Wiley-Interscience series in discrete mathematics and optimization.
	Wiley, New York, NY ; Weinheim [u.a.].
	\newblock Literaturverz. S. 721 - 747.
	
	\bibitem[Pasternak and Passy, 1972]{pasternak1972bicriterion}
	Pasternak, C. and Passy, U. (1972).
	\newblock {\em Bicriterion mathematical programs with boolean variables}.
	\newblock Technion-Israel Institute of Technology.
	
	\bibitem[Pettersson and Ozlen, 2019]{ozlen2019}
	Pettersson, W. and Ozlen, M. (2019).
	\newblock Multi-objective mixed integer programming: An objective space
	algorithm.
	\newblock {\em AIP Conference Proceedings}, 2070(1):020039.
	
	\bibitem[Przybylski et~al., 2008]{PryEtAl2008}
	Przybylski, A., Gandibleux, X., and Ehrgott, M. (2008).
	\newblock Two phase algorithms for the bi-objective assignment problem.
	\newblock {\em European Journal of Operational Research}, 185(2):509-533.
	
	\bibitem[Przybylski et~al., 2010a]{Przy10}
	Przybylski, A., Gandibleux, X., and Ehrgott, M. (2010a).
	\newblock A recursive algorithm for finding all nondominated extreme points in
	the outcome set of a multiobjective integer programme.
	\newblock {\em INFORMS Journal on Computing}, 22:371-386.
	
	\bibitem[Przybylski et~al., 2010b]{PRZYBYLSKI2010149}
	Przybylski, A., Gandibleux, X., and Ehrgott, M. (2010b).
	\newblock A two phase method for multi-objective integer programming and its
	application to the assignment problem with three objectives.
	\newblock {\em Discrete Optimization}, 7(3):149-165.
	
	\bibitem[Raith and Ehrgott, 2009]{raith09}
	Raith, A. and Ehrgott, M. (2009).
	\newblock A two-phase algorithm for the biobjective integer minimum cost flow
	problem.
	\newblock {\em http://www.esc.auckland.ac.nz/research/tech/esc-tr-661.pdf}, 36.
	
	\bibitem[Raith and Sede\~no-Noda, 2017]{raith17}
	Raith, A. and Sede\~no-Noda, A. (2017).
	\newblock Finding extreme supported solutions of biobjective network flow
	problems: An enhanced parametric programming approach.
	\newblock {\em Computers \& Operations Research}, 82.
	
	\bibitem[Sayin, 2024]{Serpil2024}
	Sayin, S. (2024).
	\newblock Supported nondominated points as a representation of the nondominated
	set: An empirical analysis.
	\newblock {\em Journal of Multi-Criteria Decision Analysis}, 31(1-2):e1829.
	
	\bibitem[Schulze, 2017]{schulze2017new}
	Schulze, B. (2017).
	\newblock {\em New perspectives on multi-objective knapsack problems}.
	\newblock PhD thesis, Dissertation, Wuppertal, Universit\"at Wuppertal, 2017.
	
	\bibitem[Schulze et~al., 2019]{Schulze2019}
	Schulze, B., Klamroth, K., and Stiglmayr, M. (2019).
	\newblock Multi-objective unconstrained combinatorial optimization: a
	polynomial bound on the number of extreme supported solutions.
	\newblock {\em Journal of Global Optimization}, 74(3):495-522.
	
	\bibitem[Sede\~no-Noda and Raith, 2015]{SEDENONODA15}
	Sede\~no-Noda, A. and Raith, A. (2015).
	\newblock A dijkstra-like method computing all extreme supported non-dominated
	solutions of the biobjective shortest path problem.
	\newblock {\em Computers \& Operations Research}, 57:83-94.
	
	\bibitem[Sourd et~al., 2006]{sourd2006multi}
	Sourd, F., Spanjaard, O., and Perny, P. (2006).
	\newblock Multi-objective branch and bound. application to the bi-objective
	spanning tree problem.
	\newblock In {\em 7th international conference in multi-objective programming
		and goal programming}.
	
	\bibitem[Steuer, 1986]{SteuerBook}
	Steuer, R.~E. (1986).
	\newblock {\em Multiple criteria optimization : theory, computation, and
		application}.
	\newblock Wiley series in probability and mathematical statistics ; Applied
	probability and mathematical statistics. Wiley, New York ;.
	
	\bibitem[{The Luc}, 1995]{dinh1995}
	{The Luc}, D. (1995).
	\newblock On the properly efficient points of nonconvex sets.
	\newblock {\em European Journal of Operational Research}, 86(2):332-336.
	
	\bibitem[Tuyttens et~al., 2000]{tuyttens2000performance}
	Tuyttens, D., Teghem, J., Fortemps, P., and Nieuwenhuyze, K.~V. (2000).
	\newblock Performance of the mosa method for the bicriteria assignment problem.
	\newblock {\em Journal of Heuristics}, 6:295-310.
	
	\bibitem[Vis\~ee et~al., 1998]{Visee1998}
	Vis\~ee, M., Teghem, J., Pirlot, M., and Ulungu, E. (1998).
	\newblock Two-phases method and branch and bound procedures to solve the
	bi-objective knapsack problem.
	\newblock {\em Journal of Global Optimization}, 12(2):139-155.
	
	\bibitem[\"Ozpeynirci and K\"oksalan, 2010]{oezpeynirci10exact}
	\"Ozpeynirci, Ã. and K\"oksalan, M. (2010).
	\newblock An exact algorithm for finding extreme supported nondominated points
	of multiobjective mixed integer programs.
	\newblock {\em Manage. Sci.}, 56(12):2302-2315.
	
\end{thebibliography}
	
%\end{thebibliography}

\end{document}